\input amstex
\input amssym
\documentstyle{amsppt}
\magnification=1200 \magnification=\magstep1 \hsize=6truein
\vsize=8.4truein \baselineskip 0.7cm
\parindent=0.8cm

\NoBlackBoxes \TagsOnLeft \magnification=\magstep1

\centerline {\bf Noether's Problem for Some $p$-groups}

\document
\baselineskip 0.7cm
\parindent=0.8cm

\bigskip

\noindent $\mskip 120mu$ Shou-Jen Hu $\mskip 180mu$  Ming-chang Kang $\mskip 20mu$

\noindent {$\mskip 80mu$ Department of Mathematics \ \  and \ \ Department
of Mathematics}

\noindent $\mskip 100mu$ Tamkang University $\mskip 110mu$ National Taiwan
University

\noindent $\mskip 120mu$Tamsui, Taiwan $\mskip 180mu$ Taipei, Taiwan

\noindent {$\mskip 370mu$ E-mail:kang\@math.ntu.edu.tw}

\vskip 20mm

\noindent {\bf Abstract}  Let $K$ be any field and $G$ be a finite group. Let $G$ act
on the rational function field $K(x_g: \ g \in G)$ by $K$-automorphisms defined by $g
\cdot x_h=x_{gh}$ for any $g, \ h \in G$. Noether's problem asks whether the fixed
field $K(G)= K(x_g: \ g \in G)^G$ is rational (=purely transcendental) over $K$.  We
will prove that if $G$ is a non-abelian $p$-group of order $p^n$ containing a cyclic
subgroup of index $p$ and $K$ is any field containing a primitive $p^{n-2}$-th root of
unity, then $K(G)$ is rational over $K$.  As a corollary, if $G$ is a non-abelian
$p$-group of order $p^3$ and $K$ is a field containing a primitive $p$-th root of
unity, then $K(G)$ is rational.

\vskip 3cm

\noindent $\underline {\mskip 200mu}$



\newpage

\noindent {\bf \S 1. Introduction}

Let $K$ be any field and $G$ be a finite group.  Let $G$ act on the rational function
field $K(x_g:g \in G)$ by $K$-automorphisms such that $g \cdot x_h=x_{gh}$ for any $g,
\ h \in G$. Denote by $K(G)$ the fixed field $K(x_g:g \in G)^G$. Noether's problem asks
whether $K(G)$ is rational (=purely transcendental) over $K$.  Noether's problem for
abelian groups was studied by Swan, Voskresenskii, Endo,  Miyata and Lenstra, etc.  See
the survey article \cite {Sw} for more details. Consequently we will restrict our
attention to the non-ableian case in this article.

First we will recall several results of Noether's problem for non-abelian
$p$-groups.

\proclaim {Theorem 1.1} {\rm (Chu and Kang \cite {CK, Theorem 1.6})} Let $G$ be a
non-abelian $p$-group of order $\le p^4$ and exponent $p^e$. Assume that $K$ is any
field such that either {\rm (i)} char $K=p >0$, or {\rm (ii)} char $K \ne p$ and $K$
contains a primitive $p^e$-th root of unity.  Then $K(G)$ is rational over $K$.
\endproclaim

\proclaim {Theorem 1.2} {\rm (\cite {Ka2, Theorem 1.5})}  Let $G$ be a non-abelian
metacyclic $p$-group of exponent $p^e$.  Assume that $K$ is any field such that either
{\rm (i)} char $K=p
>0$, or {\rm (ii)} char $K \ne p$ and $K$ contains a primitive $p^e$-th root of unity.  The $K(G)$
is rational over $K$. \endproclaim

\proclaim {Theorem 1.3} {\rm (Saltman \cite {Sa1})}  Let $K$ be any field with char $K
\ne p$ {\rm ( }in particular, $K$ may be any algebraically closed field with char $K
\neq p$ {\rm )}. There exists a non-abelian $p$-group $G$ of order $p^9$ such that
$K(G)$ is not rational over $K$.
\endproclaim

\proclaim {Theorem 1.4} {\rm (Bogomolov \cite {Bo})}  There exists a non-abelian
$p$-group $G$ of order $p^6$ such that $\Bbb C(G)$ is not rational over $\Bbb C$.
\endproclaim

All the above theorems deal with fields $K$ containing enough roots of
unity.  For a field $K$ which doesn't have enough roots of unity, so far as
we know, the only two known cases are the following Theorem 1.5 and Theorem
1.6.

\proclaim {Theorem 1.5} {\rm (Saltman \cite {Sa2, Theorem 1})}  Let $G$ be
a non-abelian $p$-group of order $p^3$.  Assume that $K$ is any field such
that either {\rm (i)} char $K =p>0$ or {\rm (ii)} char $K \ne p$ and $K$
contains a primitive $p$-th root of unity.  Then $K(G)$ is stably rational
over $K$. \endproclaim

\proclaim {Theorem 1.6} {\rm (Chu, Hu and Kang \cite {CHK; Ka1})}  Let $K$ be any
field.  Suppose that $G$ is a non-abelian group of order 8 or 16. Then $K(G)$ is
rational over $K$ except when $G=Q$, the generalized quaternion group of order 16 {\rm
(}see Theorem 1.9 for its definition {\rm)}.  When $G=Q$ and $K(\zeta)$ is cyclic over
$K$ where $\zeta$ is an primitive 8-th root of unity, then $K(G)$ is also rational over
$K$.
\endproclaim

We will remark that, if $G=Q$ is the generalized quaternion group of order 16, then
$\Bbb Q(G)$ is not rational over $\Bbb Q$ by a theorem of Serre \cite {GMS, Theorem
34.7, p.92}. The main result of this article is the following.

 \proclaim {Theorem 1.7}  Let $G$ be a non-abelian
$p$-group of order $p^n$ such that $G$ contains a cyclic subgroup of index $p$.  Assume that $K$ is
any field such that either {\rm (i)} char $K =p>0$ or {\rm (ii)} char $K \ne p$ and
$[K(\zeta):K]=1$ or $p$ where $\zeta$ is a primitive $p^{n-1}$-th root of unity.  Then $K(G)$ is
rational over $K$.
\endproclaim

As a corollary of Theorem 1.1 and Theorem 1.7, we have \proclaim {Theorem
1.8}  Let $G$ be a non-ableian $p$-group of order $p^3$.  Assume that $K$
is any field such that either {\rm (i)} char $K =p>0$ or {\rm (ii)} char $K
\ne p$ and $K$ contains a primitive $p$-th root of unity.  Then $K(G)$ is
rational over $K$. \endproclaim

Noether's problem is studied for the inverse Galois problem and the construction of a
generic Galois $G$-extension over $K$.  See \cite {DM} for details.

We will describe the main ideas of the proof of Theorem 1.7 and Theorem 1.8.  All the $p$-groups
containing cyclic subgroups of index $p$ are classified by the following theorem.

\proclaim {Theorem 1.9} {\rm (\cite {Su, p.107})}  Let $G$ be a non-ableian
$p$-group of order $p^n$ containing a cyclic subgroup of index $p$. \roster
\item "(i)" If $p$ is an odd prime number, then $G$ is isomorphic to
$M(p^n)$; and \item "(ii)"  If $p=2$, then $G$ is isomorphic to $M(2^n)$,
$D(2^{n-1}), \ SD(2^{n-1})$ where $n \ge 4$, and $\ Q(2^n)$ where $n \ge 3$
\endroster such that $$\aligned M & (p^n)=<\sigma, \tau: \ \sigma^{p^{n-1}}=\tau^p=1, \
\tau^{-1} \sigma \tau=\sigma^{1+p^{n-2}} >,\\ D & (2^{n-1})=<\sigma, \tau: \
\sigma^{2^{n-1}}=\tau^2=1, \ \tau^{-1} \sigma \tau=\sigma^{-1} >,\\ S &
D(2^{n-1})=<\sigma, \tau: \ \sigma^{2^{n-1}}=\tau^2=1, \ \tau^{-1} \sigma
\tau=\sigma^{-1+2^{n-2}}>,\\ Q & (2^n)=<\sigma, \tau: \ \sigma^{2^{n-1}}=\tau^4=1, \
\sigma ^{2^{n-2}}=\tau^2, \ \tau^{-1} \sigma \tau=\sigma^{-1}>. \endaligned$$
\endproclaim

The groups $M(p^n), \ D(2^{n-1}), \ SD(2^{n-1}), \ Q(2^n)$ are called the modular
group, the dihedral  group, the quasi-dihedral group and the generalized quaternion
group respectively.

Thus we will concentrate on the rationality of $K(G)$ for $G=M(p^n), \ D(2^{n-1}),$
$SD(2^{n-1}), \ Q(2^n)$ with the assumption that $[K(\zeta) : K] =1 $ or $p$ where $G$
is a group of exponent $p^e$ and $\zeta$ is a primitive $p^e$-th root of unity.  If
$\zeta \in K$, then Theorem 1.7 follows from Theorem 1.2.  Hence we may assume that
$[K(\zeta):K]=p$.  If $p$ is an odd prime number, the condition on $[K(\zeta):K]$
implies that $K$ contains a primitive $p^{e-1}$-th root of unity.  If $p=2$, the
condition $[K(\zeta):K]=2$ implies that $\lambda(\zeta)=- \zeta, \ \pm \zeta^{-1}$
where $\lambda$ is a generator of the Galois group of $K(\zeta)$ over $K$. (The case
$\lambda (\zeta)= - \zeta$ is equivalent to that the primitive $2^{e-1}$-th root of
unity belongs to $K$.)  In case $K$ contains a primitive $p^{e-1}$-th root of unity, we
construct a faithful representation $G \longrightarrow GL(V)$ such that dim $V=p^2$ and
$K(V)$ is rational over $K$.  For the remaining cases i.e. $p=2$, we will add the root
$\zeta$ to the ground field $K$ and show that $K(G)=K(\zeta)(G)^{<\lambda>}$ is
rational over $K$. In the case $p=2$ we will construct various faithful representations
according to the group $G=M(2^n), \ D(2^{n-1}), \ SD(2^{n-1}), \ Q(2^n)$ and the
possible image $\lambda (\zeta)$ because it seems that a straightforward imitation of
the case for $K$ containing a primitive $p^{e-1}$-th root of unity doesn't work.

We organize this article as follows.  Section 2 contains some preliminaries which will
be used subsequently.  In Section 3, we first prove Theorem 1.7 for the case when $K$
contains a primitive $p^{e-1}$-th root of unity.  This result will be applied to prove
Theorem 1.8.  In Section 4 we continue to complete the proof of Theorem 1.7.  The case
when char $K =p>0$ will be taken care by the following theorem due to Kuniyoshi.

\proclaim {Theorem 1.10} {\rm (Kuniyoshi \cite {CK, Theorem 1.7})}  If char $K=p >0$ and $G$ is a
finite $p$-group, then $K(G)$ is rational over $K$.
\endproclaim

Standing Notations.  The exponent of a finite group, denoted by exp$(G)$, is defined as
exp$(G)= \max \{{\text {ord}}(g): g \in G\}$ where ord$(g)$ is the order of the element
$g$.  Recall the definitions of modular groups, dihedral groups, quasi-dihedral groups
and generalized quaternian groups which are defined in Theorem 1.9.  If $K$ is a field
with char $K=0$ or char $K \nmid m$, then $\zeta _m$ denotes a primitive $m$-th root of
unity in some extension field of $K$.  If $L$ is any field and we write $L(x,y)$, $L(x,
y, z)$ without any explanation, we mean that these fields $L(x, y)$, $L(x, y, z)$ are
rational function fields over $K$.

\

\noindent {\bf \S 2. Generalities}

We list several results which will be used in the sequel.

\proclaim {Theorem 2.1} {\rm ([CK, Theorem 4.1])}  Let $G$ be a finite
group acting on $L(x_1, \cdots, x_m)$, the rational function field of $m$
variables over a field $L$ such that \roster \item "(i)" for any $\sigma
\in G$, $\sigma (L) \subset L$;
\item "(ii)" the restriction of the action of $G$ to $L$ is faithful;

\item "(iii)" for any $\sigma \in G$,
$$\pmatrix \sigma (x_1) \\ \vdots \\ \sigma (x_m)\endpmatrix = A(\sigma) \pmatrix x_1 \\ \vdots \\
x_m \endpmatrix +B(\sigma)$$ \endroster where $A(\sigma) \in GL_m(L)$ and $B(\sigma)$ is an $m
\times 1$ matrix over $L$. Then there exist $z_1, \cdots, z_m \in L(x_1, \cdots, x_m)$ so that
$L(x_1 \cdots, x_m)=L(z_1, \cdots, z_m)$ with $\sigma(z_i)=z_i$ for any $\sigma \in G$, any $1 \le
i \le m$.  \endproclaim

\proclaim {Theorem 2.2} {\rm ([AHK, Theorem 3.1])} Let $G$ be a finite group acting on
$L(x)$, the rational function field of one variable over a field $L$.  Assume that, for
any $\sigma \in G$, $\sigma(L) \subset L$ and $\sigma (x)=a_{\sigma} x+ b_{\sigma}$ for
any $a_{\sigma}, \ b_{\sigma} \in L$ with $a_{\sigma} \ne 0$.  Then $L(x)^G=L^G(z)$ fr
some $z \in L[x]$. \endproclaim

\proclaim {Theorem 2.3} {\rm ([CHK, Theorem 2.3])}  Let $K$ be any field, $K(x,y)$ the
rational function field of two variables over $K$, and $a, \ b \in K \setminus \{0\}$.
If $\sigma$ is a $K$-automorphism on $K(x,y)$ defined by $\sigma (x)=a/x$, $\sigma
(y)=b/y$, then $K(x, y)^{<\sigma>}= K(u,v)$ where $$\displaystyle u = \frac {x-\dfrac
ax}{xy-\dfrac {ab}{xy}},\quad v=\frac {y- \dfrac by}{xy-\dfrac {ab}{xy}}.$$

Moreover, $x+(a/x)=(-bu^2+av^2+1)/v$, $y+(b/y)=(bu^2-av^2+1)/u$,
$xy+(ab/(xy))=(-bu^2-av^2+1)/(uv)$.\endproclaim

\proclaim {Lemma 2.4}  Let $K$ be any field whose prime field is denoted by $\Bbb F$.
Let $m \ge 3$ be an integer.  Assume that char $\Bbb F \ne 2$, $[K(\zeta_{2^m}) : K]
=2$ and $\lambda(\zeta_{2^m})=\zeta_{2^m}^{-1} ($resp.
$\lambda(\zeta_{2^m})=-\zeta_{2^m}^{-1})$ where $\lambda$ is the non-trivial
$K$-automorphism on $K(\zeta_{2^m})$.  Then $K(\zeta_{2^m})=K(\zeta_4)$ and $K \bigcap
\Bbb F(\zeta_4)= \Bbb F$.
\endproclaim

\demo {Proof}  Since $m \ge 3$, it follows that $\lambda(\zeta _4)=\zeta_4^{-1}$ no matter whether
$\lambda(\zeta_{2^m})=\zeta_{2^m}^{-1}$ or $- \zeta_{2^m}^{-1}$.  Hence $\lambda (\zeta_4) \ne
\zeta_4$.  It follows that $\zeta_4 \in K(\zeta_{2^m}) \setminus K$.  Thus $K(\zeta_{2^m})=K(\zeta
_4)$.  In particular, $\zeta_4 \notin \Bbb F$. Since $[K(\zeta_4):K]=2$ and $[\Bbb F(\zeta_4): \Bbb
F]=2$, it follows that $K \bigcap \Bbb F(\zeta_4)= \Bbb F$. \qed \enddemo

\

\noindent {\bf \S 3. Proof of Theorem 1.8}

Because of Theorem 1.10 we will assume that char $K \ne p$ for any field
$K$ considered in this section.

\proclaim {Theorem 3.1}  Let $p$ be any prime number, $G=M(p^n)$ the modular group of order $p^n$
where $n \ge 3$ and $K$ be any field containing a primitive $p^{n-2}$-th root of unity.  Then
$K(G)$ is rational over $K$.
\endproclaim

\demo {Proof}  Let $\xi$ be a primitive $p^{n-2}$-th root of unity in $K$.

Step 1.

Let $\bigoplus _{g \in G} K \cdot x(g)$ be the representation space of the
regular representation of $G$.

Define $$v= \sum_{0 \le i \le p^{n-2} -1} \xi^{-i} [ x(\sigma^{ip})+x(\sigma^{ip} \tau)
+ \cdots +x(\sigma ^{ip} \tau ^{p-1})].$$

Then $\sigma ^p(v)= \xi v$ and $\tau (v)=v$.

Define $x_i=\sigma ^i v$ for $0 \le i \le p-1$.  We note that $\sigma: x_0 \mapsto x_1 \mapsto
\cdots \mapsto x_{p-1} \mapsto \xi x_0$ and $\tau: x_i \mapsto \eta^{-i} x_i$ where $\eta =
\xi^{p^{n-3}}$.

Applying Theorem 2.1 we find that, if $K(x_0, x_1, \cdots, x_{p-1})^G$ is rational over
$K$, then $K(G)=K(x(g): g \in G)^G$ is also rational over $K$.

Step 2.

Define $y_i=x_i/x_{i-1}$ for $1 \le i \le p-1$.  Then $K(x_0, x_1, \cdots,
x_{p-1})=K(x_0, y_1, $ $\cdots, y_{p-1})$ and $\sigma: x_0 \mapsto y_1x_0, y_1 \mapsto
y_2 \mapsto \cdots \mapsto y_{p-1} \mapsto \xi /(y_1 \cdots y_{p-1})$, $\tau: x_0
\mapsto x_0,$ $ \ y_i \mapsto \eta^{-1} y_i$.  By Theorem 2.2, if $K(y_1, \cdots,
y_{p-1})^G$ is rational over $K$, so is $K(x_0, y_1,$ $ \cdots, y_{p-1})^G$ over $K$.

Define $u_i=y_i/y_{i-1}$ for $2 \le i \le p-1$.  Then $K(y_1, \cdots, y_{p-1})=K(y_1,
u_2, \cdots,$ $ u_{p-1})$ and $\sigma: y_1 \mapsto y_1u_2, \ u_2 \mapsto u_3 \mapsto
\cdots \mapsto u_{p-1} \mapsto \xi/(y_1 y_2 \cdots y_{p-2}y_{p-1}^2)=\xi/$ $(y_1^p
u_2^{p-1} u_3^{p-2} \cdots u_{p-1}^2)$, $\tau: y_1 \mapsto \eta^{-1}y_1, u_i \mapsto
u_i$ for $2 \le i \le p-1$.  Thus $K(y_1,$ $ u_2, \cdots, u_{p-1})^{<\tau>}=K(y_1^p,
u_2, \cdots, u_{p-1})$.

Define $u_1=\xi^{-1}y_1^p$.  Then $\sigma: u_1 \mapsto u_1 u_2^p, u_2 \mapsto u_3
\mapsto \cdots \mapsto 1/(u_1u_2^{p-1}\cdots u_{p-1}^2)$ $\mapsto u_1
u_2^{p-2}u_3^{p-3} \cdots u_{p-2}^2u_{p-1} \mapsto u_2$.

Define $w_1=u_2$, $w_i=\sigma^{i-1}(u_2)$ for $2 \le i \le p-1$.  Then $K(u_1, u_2,
\cdots, u_{p-1})=K(w_1, w_2, \cdots, w_{p-1})$.  It follows that $K(y_1, \cdots,
y_{p-1})^G=\{K(y_1, \cdots, y_{p-1})^{<\tau>}\}^{<\sigma>}$ $=K(w_1,w_2, \cdots,
w_{p-1})^{<\sigma>}$ and $\sigma : w_1 \mapsto w_2 \mapsto \cdots \mapsto w_{p-1}
\mapsto 1/(w_1w_2 \cdots w_{p-1})$.

Step 3.

Define $T_0=1+w_1+w_1w_2+\cdots +w_1w_2 \cdots w_{p-1}$, $T_1=(1/T_0)-(1/p)$,
$T_{i+1}=(w_1w_2\cdots w_i/T_0)-(1/p)$ for $1 \le i \le p-1$.  Thus $K(w_1,
\cdots, w_{p-1})=K(T_1, \cdots, T_p)$ with $T_1 +T_2+\cdots+T_p=0$ and $\sigma
: T_1 \mapsto T_2 \mapsto \cdots \mapsto T_{p-1} \mapsto T_p \mapsto T_0$.

Define $s_i= \sum_{1 \le j \le p} \eta^{-ij}Tj$ for $1 \le i \le p-1$.  Then $K(T_1, T_2, \cdots,
T_p)=K(s_1, s_2, \cdots, s_{p-1})$ and $\sigma: s_i \mapsto \eta^is_i$.  Clearly $K(s_1, \cdots,
s_{p-1})^{<\sigma>}$ is rational over $K$. \qed \enddemo

Proof of Theorem 1.8.

If $p \ge 3$, a non-abelian $p$-group of order $p^3$ is either of exponent $p$ or contains a cyclic
subgroup of index $p$ (see \cite {CK, Theorem 2.3}).  The rationality of $K(G)$ of the first group
follows from Theorem 1.1 while that of the second group follows from the above Theorem 3.1.  If
$p=2$, the rationality of $K(G)$ is a consequence of Theorem 1.6. \qed

The method used in the proof of Theorem 3.1 can be applied to other groups, e.g. $D(2^{n-1}), \
Q(2^n), \ SD(2^{n-1})$.  The following results will be used in the proof of Theorem 1.7.

\proclaim {Theorem 3.2}  Let $G= D(2^{n-1})$ or $Q(2^n)$ with $n \ge 4$.  If $K$ is a field
containing a primitive $2^{n-2}$-th root of unity, then $K(G)$ is rational over $K$. \endproclaim

\demo {Proof}  Let $\xi$ be a primitive $2^{n-2}$-th root of unity in $K$.

Let $\bigoplus _{g \in G} K \cdot x(g)$ be the representation space of the regular representation
of $G$.

Define $$v = \sum _{0 \le i \le 2^{n-2} -1} \xi ^{-i} x(\sigma ^{2i}).$$

Then $\sigma ^2(v)= \xi v$.

Define $x_0=v, \ x_1=\sigma \cdot v, \ x_2=\tau \cdot v, \ x_3=\tau \sigma \cdot v$.  We find that
$$\aligned \sigma &: x_0 \mapsto x_1 \mapsto \xi x_0, \ x_2 \mapsto \xi^{-1}x_3, \ x_3 \mapsto x_2,
\\ \tau &: x_0 \mapsto x_2 \mapsto \epsilon x_0, \ x_1 \mapsto x_3 \mapsto \epsilon x_1
\endaligned$$
where $\epsilon =1$ if $G=D(2^{n-1})$, and $\epsilon =-1$ if $G=Q(2^n)$.

By Theorem 2.1 it suffices to show that $K(x_0, x_1, x_2, x_3)^G$ is rational over $K$.

Since $\sigma ^2(x_i)= \xi x_i$ for $i=0, 1$, $\sigma^2(x_i)=\xi^{-1}x_j$ for $j =2,
3$, it follows that $K(x_0, x_1, x_2, x_3)^{<\sigma ^2>}=K(y_o, y_1, y_2, y_3)$ where
$y_0=x_0^{2^{n-2}}, y_1=x_1/x_0, y_2=x_0x_2, y_3=x_1x_3$. The action of $\sigma$ and
$\tau$ are given by $$\aligned \sigma &: y_0 \mapsto y_0y_1^{2^{n-2}}, y_1  \mapsto
\xi/ y_1, \ y_2 \mapsto \xi^{-1}y_3, \ y_3 \mapsto \xi y_2,
\\ \tau &: y_0 \mapsto y_0^{-1}y_2^{2^{n-2}}, \ y_1 \mapsto y_1^{-1}y_2^{-1}y_3,  \ y_2 \mapsto \epsilon y_2, \ y_3
\mapsto \epsilon y_3. \endaligned$$

Define $$z_0=y_0y_1^{2^{n-3}} y_2^{-2^{n-4}}y_3^{-2^{n-4}}, \ z_1=y_1, \ z_2=y_2^{-1}y_3, \
z_3=y_2.$$

We find that $$\aligned \sigma &: z_0 \mapsto -z_0, z_1  \mapsto \xi z_1^{-1}, \ z_2 \mapsto
\xi^{2}z_2^{-1}, \ z_3 \mapsto \xi^{-1} z_2z_3,
\\ \tau &: z_0 \mapsto z_0^{-1}, \ z_1 \mapsto z_1^{-1}z_2,  \ z_2 \mapsto z_2, \ z_3
\mapsto \epsilon z_3. \endaligned$$

By Theorem 2.2 it suffices to prove that $K(z_0, z_1, z_2)^{<\sigma, \tau>}$ is
rational over $K$.

Now we will apply Theorem 2.3 to find $K(z_0, z_1, z_2)^{<\sigma>}$ with $a=1$
and $b=z_2$.  Define $$\displaystyle u=\frac {z_0-\dfrac a{z_0}}{z_0z_1-\dfrac
{ab}{z_0z_1}}, \quad v= \frac {z_1-\dfrac b{z_1}}{z_0z_1-\dfrac
{ab}{z_0z_1}}.$$

By Theorem 2.3 we find that $K(z_0, z_1, z_2)^{<\tau>} = K(u, v, z_2)$.  The actions of
$\sigma$ on $u, \ v, z_2$ are given by $$\displaystyle \aligned \sigma :& z_2 \mapsto
\xi^2z_2^{-1},\\ &u \mapsto \frac {-z_0 +\dfrac {a}{z_0}}{\xi(\dfrac {z_1}{bz_0}-\dfrac
{z_0}{z_1})}, \quad v \mapsto \frac {\xi(\dfrac 1{z_1}-\dfrac {z_1}{b})}{\xi(\dfrac
{z_1}{bz_0}-\dfrac {z_0}{z_1})}.\endaligned$$

Define $w=u/v$.  Then $\sigma (w)=bw/ \xi = z_2w/\xi$.

Note that $$\displaystyle \sigma (u) = \frac {-z_0 +\dfrac {a}{z_0}}{\xi(\dfrac
{z_1}{bz_0}-\dfrac {z_0}{z_1})}= \frac {b}{\xi}\  \frac {z_0 - \dfrac {a}{z_0}}{\dfrac
{bz_0}{z_1}-\dfrac {az_1}{z_0}}= \frac {bu}{\xi(bu^2-av^2)}.$$

The last equality of the above formula is equivalent to the following identity

$$\displaystyle \frac {x-\dfrac ax}{\dfrac {bx}y-\dfrac {ay}x}=\frac u{bu^2-av^2}. \tag
1$$ where $x, \ y, \ u, \ v, \ a, \ b$ are the same as in Theorem 2.3.  A simple way to
verify Identity (1) goes as follows:  The right-hand side of (1) is equal to
$(y+(b/y)-(1/u))^{-1}$ by Theorem 2.3.  It is not difficult to check that the left-hand
side of (1) is equal to $(y+(b/y)-(1/u))^{-1}$.

Thus $\sigma (u)=bu/(\xi(bu^2-av^2))=z_2 u/(\xi(z_2u^2-v^2))=z_2w^2/(\xi u(z_2
w^2-1))$.

Define $T=z_2w^2/ \xi, \ X=w, \ Y=u$.  Then $K(u,v,z_2)=K(T, X, Y)$ and $\sigma : T \mapsto T, X
\mapsto A/X, \ Y \mapsto B/Y$ where $A=T, \ B=T/(\xi T-1)$.  By Theorem 2.3 it follows that $K(T,
X, Y)^{<\sigma>}$ is rational over $K(T)$.  In particular, it is rational over $K$. \qed \enddemo

\proclaim {Theorem 3.3}  Let $G= SD(2^{n-1})$ with $n \ge 4$.  If $K$ is a
field containing a primitive $2^{n-2}$-th root of unity, then $K(G)$ is
rational over $K$.\endproclaim

\demo {Proof}  The case $n=4$ is a consequence of \cite {CHK, Theorem 3.2}.
Thus we may assume $n \ge 5$ in the following proof.

The proof is quite similar to that of Theorem 3.2.

Define $v,\  x_0, \ x_1, \ x_2, \ x_3$ by the same formulae as in the proof of Theorem 3.2.  Then
$\sigma : x_0 \mapsto x_1 \mapsto \xi x_0$, $x_2 \mapsto - \xi ^{-1} x_3$, $x_3 \mapsto - x_2$,
$\tau: x_0 \mapsto x_2 \mapsto x_0$, $x_1 \mapsto x_3 \mapsto x_1$.

Define $y_0=x_0^{2^{n-2}}, \ y_1=x_1/x_0, \ y_2=x_0x_2,$ and $ \ y_3=x_1x_3$. Then
$K(x_0, x_1,$ $ x_2, x_3)^{<\sigma ^2>}=K(y_0, y_1, y_2, y_3)$ and $$\aligned \sigma &:
y_0 \mapsto y_0y_1^{2^{n-2}}, \ y_1 \mapsto \xi/y_1, \ y_2 \mapsto -\xi^{-1}y_3, y_3
\mapsto -\xi y_2, \\ \tau &: y_0 \mapsto y_0^{-1} y_2^{2^{n-2}}, \ y_1 \mapsto y_1^{-1}
y_2^{-1} y_3, \ y_2 \mapsto y_2, \ y_3 \mapsto y_3. \endaligned$$

Note that the actions of $\sigma$ and $\tau$ are the same as those in the proof
of Theorem 3.2 except for the coefficients.

Thus we may define $z_0, \ z_1, \ z_2, \ z_3$ by the same formulae as in the
proof of Theorem 3.2.

Using the assumption that $n \ge 5$, we find $$\aligned \sigma &: z_0 \mapsto -z_0, \
z_1 \mapsto \xi z_1^{-1}, \ z_2 \mapsto \xi^2z_2^{-1}, z_3 \mapsto -\xi^{-1} z_2z_3, \\
\tau &: z_0 \mapsto z_0^{-1} , \ z_1 \mapsto z_1^{-1} z_2, \ z_2 \mapsto z_2, \ z_3
\mapsto z_3.
\endaligned$$

By Theorem 2.2 it suffices to prove that $K(z_0, z_1, z_2)^{<\sigma, \tau>}$ is
rational over $K$.  But the actions of $\sigma, \ \tau$ on $z_0, \ z_1, \ z_2$
are completely the same as those in the proof of Theorem 3.2.  Hence the
result. \qed \enddemo

\

\noindent {\bf \S 4. Proof of Theorem 1.7}

We will complete the proof of Theorem 1.7 in this section.

Let $\zeta$ be a primitive $p^{n-1}$-th root of unity.  If $\zeta \in K$, then
Theorem 1.7 is a consequence of Theorem 1.2.  Thus we may assume that
$[K(\zeta):K]=p$ from now on.  Let Gal($K(\zeta)/K)=<\lambda>$ and $\lambda
(\zeta)=\zeta^a$ for some integer $a$.

If $p \ge 3$, it is easy to see that $a=1 \ (mod \ p^{n-2})$ and $\zeta ^p \in K$.  By
Theorem 1.9 the $p$-group $G$ is isomorphic to $M(p^n)$.  Apply Theorem 3.1.  We are
done.

Now we consider the case $p=2$.

By Theorem 1.9 $G$ is isomorphic to $M(2^n), \ D(2^{n-1}), \ SD(2^{n-1})$ or $Q(2^n)$.  If $G$ is a
non-abelian group of order 8, the rationality of $K(G)$ is guaranteed by Theorem 1.6.  Thus it
suffices to consider the case $G$ is a 2-group of order $\ge 16$, i.e. $n \ge 4$.

Recall that $G$ is generated by two elements $\sigma$ and $\tau$ such that $\sigma^{2^{n-1}}=1$ and
$\tau ^{-1} \sigma \tau = \sigma ^k$ where \roster
\item "(i)" $k=-1$ if $G=D(2^{n-1})$ or $Q(2^n)$,
\item "(ii)" $k=1+2^{n-2}$ if $G=M(2^n)$,
\item "(iii)" $k=-1+2^{n-2}$ if $G=SD(2^{n-1})$. \endroster

As before, let $\zeta$ be a primitive $2^{n-1}$-th root of unity and Gal$(K(\zeta)/K)= <\lambda>$
with $\lambda(\zeta)=\zeta^a$ where $a^2=1$ (mod $2^{n-1})$.  It follows that the only
possibilities of $a$ (mod $2^{n-1})$ are $a=-1$, $\pm 1+2^{n-2}$.

It follows that we have four type of groups and three choices for $\lambda(\zeta)$ and
thus we should deal with 12 situations.  Fortunately many situations behaves quite
similar. And if we abuse the terminology, we may even say that some situations are
"semi-equivariant" isomorphic (but it may not be equivariant isomorphic in the usual
sense).  Hence they obey the same formulae of changing the variables.  After every
situation is reduced to a final form we may reduce the rationality problem of a group
of order $2^n \ (n \ge 4)$ to that of a group of order 16.

Let $\bigoplus _{g \in G} K \cdot x(g)$ be the representation space of the regular
representation of $G$.  We will extend the actions of $G$ and $\lambda$ to $\bigoplus
_{g \in G}K(\zeta) \cdot x(g)$ by requiring $\rho(\zeta)=\zeta$ and $\lambda
(x(g))=x(g)$ for any $\rho \in G$.  Note that $K(G) = K(x(g): g \in G)^G= \{K(\zeta)
(x(g): g \in G)^{<\lambda>} \}^G=K(\zeta)(x(g): g\in G)^{<G, \lambda>}$.

We will find a faithful subspace $\bigoplus _{0 \le i \le 3}K(\zeta) \cdot x_i$ \ of \
$\bigoplus_{g \in G} K(\zeta) \cdot x(g)$ such that $K(\zeta)(x_0, x_1, x_2, x_3)^{<G,
\lambda>}(y_1, \cdots, y_{12})$ is rational over $K$ where each $y_i$ is fixed by $G$
and $\lambda$.  By Theorem 2.1, $K(\zeta)(x(g): g \in G)^{<G, \lambda>}=K(\zeta)(x_0,
x_1, x_2, x_3)^{<G, \lambda>}(X_1,$ $ \cdots, X_{N})$ where $N=2^n-4$ and each  $X_i$
is fixed by  $G$ and  $\lambda$.  It follows that  $K(G)$ is rational provided that
$K(\zeta)(x_0, x_1, x_2, x_3)^{<G, \lambda>}(y_1, \cdots, y_{12})$ is rational over
$K$.

Define $$v_1=\sum _{0 \le j \le 2^{n-1}-1} \zeta^{-j} x(\sigma^j), \quad v_2=\sum _{0 \le j \le
2^{n-1}-1} \zeta^{-aj} x(\sigma^j)$$ where $a$ is the integer with $\lambda (\zeta) = \zeta^a$.

We find that $\sigma : v_1 \mapsto \zeta v_1, \ v_2 \mapsto \zeta ^a v_2, \ \lambda: v_1 \mapsto
v_2 \mapsto v_1$.

Define $x_0=v_1, \ x_1 = \tau \cdot v_1, \ x_2 = v_2, \ x_3= \tau \cdot v_2$.

It follows that

$$\aligned \sigma &: x_0 \mapsto \zeta x_0, \ x_1 \mapsto \zeta^k x_1, \ x_2 \mapsto
\zeta^a x_2, \ x_3 \mapsto \zeta^{ak} x_3,\\ \lambda &: x_0 \mapsto x_2 \mapsto x_0, \
x_1 \mapsto x_3 \mapsto x_1, \ \zeta \mapsto \zeta^a,
\\ \tau &:  x_0 \mapsto x_1 \mapsto \epsilon x_0, \ x_2 \mapsto x_3 \mapsto \epsilon x_2,\\ \tau \lambda &:
x_0 \mapsto x_3 \mapsto \epsilon x_0, \ x_1 \mapsto \epsilon x_2, \ x_2 \mapsto x_1, \ \zeta
\mapsto \zeta^a\endaligned$$ where (i) $\epsilon =1$ if $G \ne Q(2^n)$, and (ii) $\epsilon =-1$ if
$G=Q(2^n)$.

Case 1. $k=-1$, i.e. $G=D(2^{n-1})$ or $Q(2^n)$.

Throughout the discussion of this case, we will adopt the convention that $\epsilon =1
$ if $G=D(2^{n-1}),$ while $\epsilon =-1$ if $G=Q(2^n)$.

Subcase 1.1.  $a=-1$, i.e. $\lambda(\zeta)= \zeta^{-1}$.

It is easy to find that $K(\zeta)(x_0, x_1, x_2,
x_3)^{<\sigma>}=K(\zeta)(x_0^{2^{n-1}}, x_0x_1, x_0x_2, x_1x_3).$

Define $$y_0=x_0^{2^{n-1}}, \ y_1=x_0x_1, \ y_2=x_0x_2, \ y_3=x_1 x_3.$$

It follows that $$\aligned \lambda &:y_0 \mapsto y_0^{-1}y_2^{2^{n-1}}, \ \ y_1 \mapsto y_1^{-1}
y_2 y_3, \ \ y_2 \mapsto y_2, \ \ y_3 \mapsto y_3, \ \ \zeta \mapsto \zeta^{-1}, \\ \tau &:y_0
\mapsto y_0^{-1} y_1 ^{2^{n-1}}, \ \ y_1 \mapsto \epsilon y_1, \ \ y_2 \mapsto y_3 \mapsto y_2.
\endaligned$$

Define  $$z_0=y_0y_1^{-2^{n-2}}y_2^{-2^{n-3}}y_3^{2^{n-3}}, \ z_1=y_2y_3, \ z_2=y_2, \ z_3=y_1 .$$

We find that $$\aligned \lambda &:z_0 \mapsto 1/z_0, \ \ z_1 \mapsto z_1, \ \ z_2
\mapsto z_2, \ \ z_3 \mapsto z_1/z_3, \ \ \zeta \mapsto \zeta^{-1}, \\ \tau &:z_0
\mapsto 1/z_0, \ \ z_1 \mapsto z_1, \ \ z_2 \mapsto z_1/z_2, \ \ z_3 \mapsto \epsilon
z_3.
\endaligned$$

It turns out the parameter $n$ does not come into play in the actions of $\lambda$ and $\tau$ on
$z_0, \ z_1, \ z_2, \ z_3$.

By Theorem 2.1 $K(G) = K(\zeta)(z_0, z_1, z_2, z_3)^{<\lambda, \tau>} (X_1, \cdots, X_N)$ where
$N=2^n-4$ and $\lambda(X_i)=\tau (X_i)=X_i$ for $1 \le i \le N$.

By Lemma 2.4 $K(\zeta)=K(\zeta_4)$ where $\lambda(\zeta_4)=\zeta_4^{-1}$.  Thus
$K(G)=K(\zeta_4)(z_0, z_1,$ $z_2, z_3)^{<\lambda, \tau>}(X_1, \cdots, X_N)$

Denote $G_4=D(8)$ or $Q(16)$.  Then $K(G_4)=K(\zeta_4)(z_0, z_1, z_2, z_3)^{<\lambda,
\tau>}(X_1, \cdots,$ $ X_{12})$.  Since $K(G_4)$ is rational over $K$ by Theorem 1.6
(see \cite {Ka1, Theorem 1.3}), it follows that $K(\zeta_4)(z_0, \cdots,
z_3)^{<\lambda, \tau>}(X_1, \cdots, X_{12})$ is rational over $K$.  Thus $K(\zeta_4)$ $
(z_0,$ $\cdots, z_3)^{<\lambda, \tau>}(X_1, \cdots, X_N)$ is rational over $K$ for
$N=2^n-4$.  The last field is nothing but $K(G)$.  Done.

Subcase 1.2.  $a=-1+2^{n-2}$, i.e. $\lambda(\zeta)=-\zeta^{-1}$.

The actions of $\sigma, \ \tau, \ \lambda, \ \tau \lambda$ are given by

$$\aligned \sigma &: x_0 \mapsto \zeta x_0, \ x_1 \mapsto \zeta^{-1} x_1, \ x_2 \mapsto
-\zeta^{-1} x_2, \ x_3 \mapsto -\zeta x_3,\\ \lambda &: x_0 \mapsto x_2 \mapsto x_0, \
x_1 \mapsto x_3 \mapsto x_1, \ \zeta \mapsto - \zeta^{-1},
\\ \tau &:  x_0 \mapsto x_1 \mapsto \epsilon x_0, \ x_2 \mapsto x_3 \mapsto \epsilon x_2,\\ \tau \lambda &:
x_0 \mapsto x_3 \mapsto \epsilon x_0, \ x_1 \mapsto \epsilon x_2, \ x_2 \mapsto x_1, \ \zeta
\mapsto -\zeta^{-1}\endaligned$$

Define $y_0=x_0^{2^{n-1}}, \ y_1 = x_0 x_1, \ y_2=x_2 x_3, \ y_3=x_0^{-1-2^{n-2}} x_3.$
Then $K(\zeta)(x_0,$ $\cdots, x_3)^{<\sigma>} = K(\zeta)(y_0, \cdots, y_3)$.  Consider
the actions of $\tau \lambda$ and $\tau$ on $K(\zeta)(y_0, \cdots,$ $ y_3)$.  We find
that $$\aligned \tau \lambda &: y_0 \mapsto y_0^{1+2^{n-2}} y_3^{2^{n-1}}, \ y_1
\mapsto \epsilon y_2 \mapsto y_1, \ y_3 \mapsto \epsilon
y_0^{-1-2^{n-3}}y_3^{-1-2^{n-2}}, \ \zeta \mapsto - \zeta ^{-1}, \\ \tau &: y_0 \mapsto
y_0^{-1} y_1^{2^{n-1}}, \ y_1 \mapsto \epsilon y_1, \ y_2 \mapsto \epsilon y_2, \ y_3
\mapsto \epsilon y_1 ^{-1-2^{n-2}} y_2 y_3^{-1}. \endaligned$$

Define $$z_0=y_1, \ z_1=y_1^{-1}y_2, \ z_2=y_0 y_1 y_2^{-1}y_3^2, \ z_3=y_0^{1+2^{n-4}} y_1
^{-2^{n-4}}y_2^{-2^{n-4}}y_3^{1+2^{n-3}}.$$

We find $$\aligned \tau \lambda &: z_0 \mapsto \epsilon z_0 z_1, \ z_1 \mapsto 1/z_1, \
z_2 \mapsto 1/z_2, \ z_3 \mapsto \epsilon z_1 ^{-1} z_2 ^{-1} z_3, \ \zeta \mapsto -
\zeta^{-1}, \\ \tau &: z_0 \mapsto \epsilon z_0, \ z_1 \mapsto z_1, \ z_2 \mapsto
1/z_2, \ z_3 \mapsto \epsilon z_1 / z_3.
\endaligned$$

By Lemma 2.4 we may replace $K(\zeta)$ in $K(\zeta)(z_0, z_1, z_2, z_3)^{<\tau \lambda, \tau>}$ by
$K(\zeta_4)$ where $\tau \lambda (\zeta_4)=\zeta_4^{-1}$.  Then we may proceed as in Subcase 1.1.
The details are omitted.

Subcase 1.3.  $a=1+2^{n-2}$, i.e. $\lambda (\zeta)= -\zeta$.

Note that $\zeta^2 \in K$ and $\zeta^2$ is a primitive $2^{n-2}$-th root of unity.  Thus we may
apply Theorem 3.2.  Done

Case 2.  $k=1+2^{n-2}$, i.e. $G=M(2^n)$.

Subcase 2.1.  $a=-1$, i.e. $\lambda (\zeta)=\zeta^{-1}$.

The actions of $\sigma, \ \tau, \ \lambda, \ \tau \lambda$ are given by  $$\aligned
\sigma &: x_0 \mapsto \zeta x_0, \ x_1 \mapsto - \zeta x_1, \ x_2 \mapsto \zeta^{-1}
x_2, \ x_3 \mapsto -\zeta^{-1} x_3,\\ \lambda &: x_0 \mapsto x_2 \mapsto x_0, \ x_1
\mapsto x_3 \mapsto x_1, \ \zeta \mapsto \zeta^{-1},
\\ \tau &:  x_0 \mapsto x_1 \mapsto x_0, \ x_2 \mapsto x_3 \mapsto  x_2,\\ \tau \lambda &:
x_0 \mapsto x_3 \mapsto x_0, \ x_1 \mapsto x_2 \mapsto x_1, \ \zeta \mapsto
\zeta^{-1}.\endaligned$$

Define $X_0=x_0, \ X_1=x_2, \ X_2=x_3, \ X_3=x_1$.  Then the actions of $\sigma, \
\tau, \ \lambda$ on $X_0, \ X_1, \ X_2, \ X_3$ are the same as those of $\sigma, \ \tau
\lambda, \ \tau,$ on $x_0, \ x_1, \ x_2, \ x_3$ in Subcase 1.2 for $D(2^{n-1})$ except
on $\zeta$.  Thus we may consider$K(\zeta)(X_0, X_1, X_2, X_3)^{<\sigma, \tau,
\lambda>}(Y_1, \cdots,$ $ Y_{12})$.  Hence the same formulae of changing the variables
in Subcase 1.2 can be copied and the same method can be used to prove that
$K(\zeta)(X_0, X_1, X_2, X_3)^{<\sigma, \tau, \lambda>}$ $(Y_1, \cdots, Y_{12})$ is
rational over $K$.

Subcase 2.2. $a=-1+2^{n-2}$, i.e. $\lambda (\zeta) = - \zeta^{-1}$.

The actions of $\sigma, \ \tau, \ \lambda, \ \tau \lambda$ are given by $$\aligned
\sigma &: x_0 \mapsto \zeta x_0, \ x_1 \mapsto -\zeta x_1, \ x_2 \mapsto - \zeta^{-1}
x_2, \ x_3 \mapsto \zeta^{-1} x_3,\\ \lambda &: x_0 \mapsto x_2 \mapsto x_0, \ x_1
\mapsto x_3 \mapsto x_1, \ \zeta \mapsto -\zeta^{-1},
\\ \tau &:  x_0 \mapsto x_1 \mapsto x_0, \ x_2 \mapsto x_3 \mapsto  x_2,\\ \tau \lambda &:
x_0 \mapsto x_3 \mapsto x_0, \ x_1 \mapsto x_2 \mapsto x_1, \ \zeta \mapsto
-\zeta^{-1}.\endaligned$$

Define $X_0=x_0, \ X_1=x_3, \ X_2=x_2, \ X_3=x_1$.  Then the actions of $\sigma, \
\tau, \ \tau \lambda$ on $X_0, \ X_1, \ X_2, \ X_3$ are the same as those of $\sigma, \
\tau \lambda, \ \tau,$ on $x_0, \ x_1, \ x_2, \ x_3$ in Subcase 1.2 for $D(2^{n-1})$.
Hence the result.

Subcase 2.3. $a=1+2^{n-2}$, i.e. $\lambda (\zeta) = - \zeta$.

Apply Theorem 3.1.

Case 3.  $k=-1+2^{n-2}$, i.e. $G=SD(2^{n-1})$.

Subcase 3.1. $a=-1$, i.e. $\lambda (\zeta) = \zeta^{-1}$.

The actions of $\sigma, \ \tau, \ \lambda, \ \tau \lambda$ are given by $$\aligned
\sigma &: x_0 \mapsto \zeta x_0, \ x_1 \mapsto -\zeta^{-1} x_1, \ x_2 \mapsto
\zeta^{-1} x_2, \ x_3 \mapsto -\zeta x_3,\\ \lambda &: x_0 \mapsto x_2 \mapsto x_0, \
x_1 \mapsto x_3 \mapsto x_1, \ \zeta \mapsto \zeta^{-1},
\\ \tau &:  x_0 \mapsto x_1 \mapsto x_0, \ x_2 \mapsto x_3 \mapsto  x_2,\\ \tau \lambda &:
x_0 \mapsto x_3 \mapsto x_0, \ x_1 \mapsto x_2 \mapsto x_1, \ \zeta \mapsto
\zeta^{-1}.\endaligned$$

Define $X_0=x_0, \ X_1=x_2, \ X_2=x_1, \ X_3=x_3$.  Then the actions of $\sigma, \ \tau
\lambda, \ \lambda $ on $X_0, \ X_1, \ X_2, \ X_3$ are the same as those of $\sigma, \
\tau \lambda, \ \tau,$ on $x_0, \ x_1, \ x_2, \ x_3$ in Subcase 1.2 for $D(2^{n-1})$
except on $\zeta$. Done.

Subcase 3.2. $a=-1+2^{n-2}$, i.e. $\lambda (\zeta) = - \zeta^{-1}$.

Define $y_0=x_0^{2^{n-1}}, \ y_1 = x_0^{1+2^{n-2}}x_1, \ y_2=x_1 ^{-1} x_2, \ y_3=x_0
^{-1} x_3$. Then $K(\zeta)(x_0, x_1, x_2, x_3)^{<\sigma>}=K(\zeta)(y_0, y_1, y_2, y_3)$
and

$$\aligned \tau &: y_0 \mapsto y_0^{-1-2^{n-2}}y_1^{2^{n-1}}, \ y_1 \mapsto
y_0^{-1-2^{n-3}} y_1 ^{1+2^{n-2}}, \ y_2 \mapsto y_3 \mapsto y_2, \\ \tau \lambda &:
y_0 \mapsto y_0 y_3^{2^{n-1}}, \ y_1 \mapsto y_1 y_2 y_3^{1+2^{n-2}}, \ y_2 \mapsto
y_2^{-1}, y_3 \mapsto y_3^{-1}, \ \zeta \mapsto -\zeta^{-1}. \endaligned$$

Define $z_0 = y_0 ^{1+2^{n-3}}y_1^{-2^{n-2}} y_2 ^{-2^{n-3}}y_3^{2^{n-3}}, \ z_1 = y_0 ^{2^{n-4}}
y_1 ^{1-2^{n-3}}y_2^{-2^{n-4}}y_3^{2^{n-4}}, \ z_2=y_2, \ z_3=y_2^{-1}y_3$.  It follows that
$K(\zeta)(y_0, y_1, y_2, y_3)=K(\zeta) (z_0, z_1, z_2, z_3)$ and $$\aligned \tau &: z_0 \mapsto
1/z_0, \ z_1 \mapsto z_1/z_0 , \ z_2 \mapsto z_2z_3, \  z_3 \mapsto 1/z_3, \\ \tau \lambda &: z_0
\mapsto z_0, \ z_1 \mapsto z_1 z_2^2 z_3, \ z_2 \mapsto 1/z_2, z_3 \mapsto 1/z_3, \ \zeta \mapsto
-\zeta^{-1}.\endaligned$$

Thus we can establish the rationality because we may replace $K(\zeta)$ by $K(\zeta_4)$ as in
Subcase 1.2.

Subcase 3.3. $a=1+2^{n-2}$, i.e. $\lambda(\zeta)= - \zeta$.

Apply Theorem 3.3.

Thus we have finished the proof of Theorem 1.7. \qed

\newpage

\centerline {\bf REFERENCES}

\roster

\noindent \item "[AHK]" H.\ Ahmad, M.\ Hajja and M.\ Kang, {\it Rationality of some projective
linear actions}, J.\ Algebra {\bf 228}(2000) 643-658.

\noindent \item "[Bo] " F.\ A.\ Bogomolov, {\it The Brauer group of quotient spaces by linear group
actions}, Math. USSR Izv. {\bf 30}(1988) 455-485.

 \noindent \item "[CHK]" H.\ Chu, S. \ J.\ Hu and  M.\ Kang, {\it
Nother's problem for dihedral 2-groups}, Comment. Math. Helv. {\bf 79}(2004) 147-159.

\noindent \item "[CK]  " H.\ Chu and M.\ Kang,  {\it Rationality of p-group actions},
J.\ Algebra {\bf 237}(2001) 673-690.

\noindent \item "[DM]" F. DeMeyer and T. McKenzie, {\it On generic polynomials}, J.
Alg. {\bf 261} (2003) 327-333.

\noindent \item "[GMS]" S. Garibaldi, A. Merkurjev and J. P. Serre, Cohomological
invariants in Galois cohomology, AMS Univ. Lecture Series vol. 28, Amer. Math. Soc.,
Providence, 2003.

\noindent \item "[Ka1] " M.\ Kang,  {\it Noether's problem for dihedral 2-groups II},
to appear in "Pacific J. Math.".

\noindent \item "[Ka2] " M.\ Kang,  {\it Noether's problem for metacyclic $p$-groups},
to appear in "Advances in Math.".

\item "[Sa1]" D. J. Saltman, {\it Noether's problem over an algebraically closed
field}, Invent. Math. {\bf 77} (1984) 71-84.

\noindent \item "[Sa2]" D.\ J. \ Saltman, {\it Galois groups of order $p^3$}, Comm. Alg.\ {\bf 15}
(1987) 1365-1373.

\noindent \item "[Su]" M. Suzuki, Group theory II, Grund. math. Wiss. Vol. 248, Springer-Verlag,
Berlin, 1986.

\noindent \item "[Sw] " R.\ G.\ Swan, {\it  Noether's problem in Galois theory}, in
"Emmy Noether in Bryn Mawr", edited by B. Srinivasan and J. Sally, Springer-Verlag,
Berlin, 1983.

\endroster

\end